\documentclass[12pt,oneside,english]{amsart}
\usepackage[T1]{fontenc}
\usepackage[latin1]{inputenc}
\usepackage{geometry}
\usepackage{setspace}
\usepackage{amssymb}

\makeatletter
 \theoremstyle{plain}
\newtheorem{thm}{Theorem}[section]
  \theoremstyle{plain}
  \newtheorem{cor}[thm]{Corollary}
  \theoremstyle{plain}
  \newtheorem{lem}[thm]{Lemma}

\numberwithin{equation}{section}

\usepackage{babel}
\makeatother
\begin{document}

\title{on freely indecomposable measures }

\author{hari bercovici and jiun-chau wang}

\date{March 5th, 2007; Revised on October 5th, 2007}

\begin{abstract}
We show that a probability measure is not a nontrivial free convolution
if it puts no mass in an interval whose endpoints are atoms. The proof
uses analytic subordination. 
\end{abstract}

\thanks{The first author was supported in part by a grant from the National
Science Foundation.}

\maketitle

\section{Introduction}

Given two probability measures $\mu,\nu$ on the real line $\mathbb{R}$,
we denote by $\mu\boxplus\nu$ their free convolution (see \cite{VDN}
for the definition of free convolution). If $\nu$ is a point mass,
then the measure $\mu\boxplus\nu$ is just a translation of $\mu$.
A measure of the form $\mu\boxplus\nu$, where neither $\mu$ nor
$\nu$ is a point mass, is said to be \emph{freely decomposable}.
Several classes of measures are known to be freely \emph{indecomposable}.
For instance, Belinschi proved in \cite{Serban1,Serban2} that a measure
with nontrivial continuous singular part is freely indecomposable.
More recently, Chistyakov and G\"{o}tze observed in \cite{CG1} that
measures with finite support are freely indecomposable (this result
also follows from the description given in \cite{BV} of the atoms
of a free convolution.) Both of these classes of measures are weakly
dense in the set of all Borel probability measures on $\mathbb{R}$. 

In this note we will prove that $\mu$ is freely indecomposable if
there are points $\alpha<\beta$ such that $\mu(\{\alpha\})>0$, $\mu(\{\beta\})>0$,
and $\mu((\alpha,\beta))=0$. We also prove analogous results for
free multiplicative convolutions $\boxtimes$ of measures defined
on the positive half-line $\mathbb{R}_{+}=[0,+\infty)$, and on the
circle $\mathbb{T}=\{\zeta\in\mathbb{C}:\,\left|\zeta\right|=1\}$.

\section{Additive Free Convolution}

Given a probability measure $\mu$ on $\mathbb{R}$, we define the
analytic function $G_{\mu}$ on $\mathbb{C}^{+}=\{ z\in\mathbb{C}:\,\Im z>0\}$
by \[
G_{\mu}(z)=\int_{-\infty}^{\infty}\frac{1}{z-t}\, d\mu(t),\qquad\Im z>0.\]
Note that the measure $\mu$ is completely determined by the imaginary
part of $G_{\mu}$. Set $\mathbb{C}^{-}=-\mathbb{C}^{+}$. A free
convolution $\mu_{1}\boxplus\mu_{2}$ is characterized analytically
by the identity \begin{equation}
G_{\mu_{1}\boxplus\mu_{2}}^{-1}(w)=G_{\mu_{1}}^{-1}(w)+G_{\mu_{2}}^{-1}(w)-\frac{1}{w},\label{eq:2.1}\end{equation}
where $G_{\mu}^{-1}$ denote the inverse of $G_{\mu}$ relative to
composition, and $w$ belongs to an appropriate Stolz angle at zero
in $\mathbb{C}^{-}$, say $\left|\Re w\right|<-\Im w<\varepsilon$
for some $\varepsilon>0$. 

It was shown by Biane \cite{Biane} (cf. also \cite{Free entropy}
for an earlier partial result) that, given measures $\mu_{1}$ and
$\mu_{2}$, there exist analytic functions $\omega_{1},\omega_{2}:\,\mathbb{C}^{+}\to\mathbb{C}^{+}$
such that \begin{equation}
G_{\mu_{1}\boxplus\mu_{2}}(z)=G_{\mu_{1}}(\omega_{1}(z))=G_{\mu_{2}}(\omega_{2}(z)),\qquad z\in\mathbb{C}^{+}.\label{eq:2.2}\end{equation}
The functions $\omega_{1},\omega_{2}$ are uniquely determined, and
they satisfy\[
\lim_{y\rightarrow+\infty}\frac{\omega_{j}(iy)}{iy}=1,\qquad j=1,2.\]
Moreover, as observed in \cite{BV}, relation (2.1) can be rewritten
as \begin{equation}
\omega_{1}(z)+\omega_{2}(z)=z+\frac{1}{G_{\mu_{1}\boxplus\mu_{2}}(z)},\qquad z\in\mathbb{C}^{+}.\label{eq:2.3}\end{equation}

The following result is proved in \cite{BV}.

\begin{thm}
Assume $\alpha$ is an atom of the measure $\mu_{1}\boxplus\mu_{2}$.
Then
\begin{enumerate}
\item the limits $\alpha_{j}=\lim_{\varepsilon\downarrow0}\omega_{j}(\alpha+i\varepsilon)$
exist, $j=1,2$.
\item $\alpha_{1}+\alpha_{2}=\alpha$.
\item $\mu_{1}(\{\alpha_{1}\})+\mu_{2}(\{\alpha_{2}\})=(\mu_{1}\boxplus\mu_{2})(\{\alpha\})+1$.
\item \[
\lim_{\varepsilon\downarrow0}\frac{\omega_{j}(\alpha+i\varepsilon)-\alpha_{j}}{i\varepsilon}=\frac{\mu_{j}(\{\alpha_{j}\})}{(\mu_{1}\boxplus\mu_{2})(\{\alpha\})},\; j=1,2.\]

\end{enumerate}
\end{thm}
Part (1) actually occurs in the proof of Theorem 7.4 of \cite{BV},
while (4) is only implicit in that proof. The relevant calculation
goes as follows for $j=1$:\[
\frac{\omega_{1}(\alpha+i\varepsilon)-\alpha_{1}}{i\varepsilon}=\frac{(\omega_{1}(\alpha+i\varepsilon)-\alpha_{1})G_{\mu_{1}}(\omega_{1}(\alpha+i\varepsilon))}{(\alpha+i\varepsilon-\alpha)G_{\mu_{1}\boxplus\mu_{2}}(\alpha+i\varepsilon)}.\]
By Lemma 7.1 in \cite{BV}, the numerator and denominator of the last
fraction converge respectively to $\mu_{1}(\{\alpha_{1}\})$ and $(\mu_{1}\boxplus\mu_{2})(\{\alpha\})$
as $\varepsilon\rightarrow0^{+}$. 

\begin{cor}
Assume that $\alpha$ and $\beta$ are atoms of $\mu_{1}\boxplus\mu_{2}$,
and write them as \[
\alpha=\alpha_{1}+\alpha_{2},\:\beta=\beta_{1}+\beta_{2}\]
as in the preceding theorem. Then either $\alpha_{1}=\beta_{1}$ or
$\alpha_{2}=\beta_{2}$. 
\end{cor}
\begin{proof}
If $\alpha_{1}\neq\beta_{1}$ and $\alpha_{2}\neq\beta_{2}$, then\begin{gather*}
2<2+(\mu_{1}\boxplus\mu_{2})(\{\alpha\})+(\mu_{1}\boxplus\mu_{2})(\{\beta\})=\qquad\qquad\qquad\quad\\
\mu_{1}(\{\alpha_{1}\})+\mu_{1}(\{\beta_{1}\})+\mu_{2}(\{\alpha_{2}\})+\mu_{2}(\{\beta_{2}\})\leq2,\end{gather*}
a contradiction. 
\end{proof}
From this point on, we will assume that $\mu_{1}$ and $\mu_{2}$
are not point masses, $\mu_{1}\boxplus\mu_{2}$ has two atoms $\alpha<\beta$,
and $(\mu_{1}\boxplus\mu_{2})((\alpha,\beta))=0$. Let us write $\alpha=\alpha_{1}+\alpha_{2}$
and $\beta=\beta_{1}+\beta_{2}$ as in Corollary 2.2. Exchanging $\mu_{1}$
and $\mu_{2}$ if necessary, we may assume that $\alpha_{2}=\beta_{2}$.
Furthermore, replacing $\mu_{2}$ by $\mu_{2}\boxplus\delta_{-\alpha_{2}}$,
we may assume that $\alpha_{2}=\beta_{2}=0$ so that $\alpha=\alpha_{1}$
and $\beta=\beta_{1}$ are atoms of $\mu_{1}$ as well as $\mu_{1}\boxplus\mu_{2}$. 

\begin{lem}
The functions $\omega_{1}$ and $\omega_{2}$ can be extended meromorphically
across $(\alpha,\beta)$. Both continuations are real-valued on $(\alpha,\beta)$,
with the exception of at most one pole.
\end{lem}
\begin{proof}
Since $\omega_{1},\omega_{2}$, and $1/G_{\mu_{1}\boxplus\mu_{2}}$
take values in $\mathbb{C}^{+}$, they have Nevanlinna representations
\[
\omega_{j}(z)=r_{j}+z+\int_{-\infty}^{\infty}\frac{1+tz}{t-z}\, d\sigma_{j}(t),\qquad j=1,2,\]
\[
\frac{1}{G_{\mu_{1}\boxplus\mu_{2}}(z)}=s+z+\int_{-\infty}^{\infty}\frac{1+tz}{t-z}\, d\sigma(t),\]
for $z\in\mathbb{C}^{+}$, where $r_{j},s\in\mathbb{R}$, and $\sigma_{j},\sigma$
are finite positive Borel measures on $\mathbb{R}$. The identity
(2.3) implies that $\sigma_{1}+\sigma_{2}=\sigma$. Now, the assumption
that $(\mu_{1}\boxplus\mu_{2})((\alpha,\beta))=0$ implies that $G_{\mu_{1}\boxplus\mu_{2}}$
can be continued analytically across $(\alpha,\beta)$, and this continuation,
which we still denote by $G_{\mu_{1}\boxplus\mu_{2}}$, is real-valued
and strictly decreasing on $(\alpha,\beta)$. 

Since $\alpha,\beta$ are atoms of $\mu_{1}\boxplus\mu_{2}$, we have
\[
\lim_{t\downarrow\alpha}G_{\mu_{1}\boxplus\mu_{2}}(t)=+\infty,\]
and\[
\lim_{t\uparrow\beta}G_{\mu_{1}\boxplus\mu_{2}}(t)=-\infty.\]
Hence, there exists a unique $\gamma\in(\alpha,\beta)$ so that $G_{\mu_{1}\boxplus\mu_{2}}(\gamma)=0$.
It follows that the function $1/G_{\mu_{1}\boxplus\mu_{2}}$ can be
extended meromorphically to $(\alpha,\beta)$, with a single simple
pole at $\gamma$. This means that $\sigma((\alpha,\gamma))=\sigma((\gamma,\beta))=0$
and $\sigma(\{\gamma\})>0$. Therefore $\sigma_{j}((\alpha,\gamma))=\sigma_{j}((\gamma,\beta))=0$,
and this implies the claimed properties of $\omega_{j}$. 
\end{proof}
We will need one more detail about the boundary behavior of $\omega_{j}$
which is given by the following result. 

\begin{lem}
Let $\omega:\,\mathbb{C}^{+}\to\mathbb{C}^{+}$ be an analytic function
and $\alpha,\gamma\in\mathbb{R}$. Assume that
\begin{enumerate}
\item $\omega$ can be continued analytically across $(\alpha,\gamma)$;
and also use $\omega$ to denote this continuation.
\item $\omega$ is real-valued on $(\alpha,\gamma)$.
\item $\lim_{\varepsilon\downarrow0}\omega(\alpha+i\varepsilon)=\alpha$.
\item the limit\[
a=\lim_{\varepsilon\downarrow0}\frac{\omega(\alpha+i\varepsilon)-\alpha}{i\varepsilon}\]
is finite. 
\end{enumerate}
Then \[
\lim_{z\rightarrow\alpha,\,\Re z>\alpha}\frac{\omega(z)-\alpha}{z-\alpha}=a,\]
In particular, $\lim_{t\downarrow\alpha}\omega(t)=\alpha$. 

\end{lem}
\begin{proof}
Note that the limit in (4) always exists, and belongs to $(0,+\infty]$.
This is called the \emph{Julia-Carath\'{e}odory derivative} of $\omega$
at $\alpha$ (see Exercises 6 and 7 in \cite[Chapter I]{Garnett}).
Let us also note that the function $\omega$ is strictly increasing
on $(\alpha,\gamma)$. It will be easier to work with the function
\[
\widetilde{\omega}\equiv\varphi\circ\omega\circ\varphi^{-1},\]
where \[
\varphi(z)=\frac{z-\gamma}{z-\alpha},\qquad z\in\mathbb{C}^{+}.\]
The assumptions means that 
\begin{enumerate}
\item $\widetilde{\omega}$ can be extended meromorphically to $(-\infty,0)$.
\item $\widetilde{\omega}$ is real-valued on $(-\infty,0)$, with the exception
of at most one pole, say, at $t_{0}\in(-\infty,0)$. 
\item $\lim_{y\rightarrow+\infty}\widetilde{\omega}(iy)=\infty$, and
\item the limit \[
\lim_{y\rightarrow+\infty}\frac{\widetilde{\omega}(iy)}{iy}=\frac{1}{a}\neq0.\]
 
\end{enumerate}
The Nevanlinna integral representation for $\widetilde{\omega}$ is
therefore \[
\widetilde{\omega}(z)=r+\frac{z}{a}+\frac{1+t_{0}z}{t_{0}-z}\sigma(\{ t_{0}\})+\int_{0}^{\infty}\frac{1+tz}{t-z}\, d\sigma(t),\qquad z\in\mathbb{C}^{+},\]
where $r\in\mathbb{R}$, and $\sigma$ is a finite positive Borel
measure on $[0,+\infty)$. Observe now that\[
\left|\frac{t}{t-z}\right|\leq1,\qquad t\in[0,+\infty),\:\Re z<0,\]
and the dominated convergence theorem easily yields \[
\lim_{z\rightarrow\infty,\,\Re z<0}\frac{\widetilde{\omega}(z)}{z}=\frac{1}{a}.\]
This is immediately seen to be equivalent to the conclusion of the
lemma. 
\end{proof}
We are now ready for the main result of this section. We denote by
$\delta_{t}$ the unit point mass at $t$.

\begin{thm}
Let $\mu_{1},\mu_{2}$ be probability measures on $\mathbb{R}$, and
$\alpha<\beta$. If $\alpha$ and $\beta$ are atoms of $\mu_{1}\boxplus\mu_{2}$,
and $(\mu_{1}\boxplus\mu_{2})((\alpha,\beta))=0$, then either $\mu_{1}$
or $\mu_{2}$ is a point mass. 
\end{thm}
\begin{proof}
Assume to the contrary that neither $\mu_{1}$ nor $\mu_{2}$ are
point masses. We may, and do, assume that $\alpha$ and $\beta$ are
atoms of the measure $\mu_{1}$. With the notation used earlier, Lemma
2.4 and (2.3) imply that \[
\lim_{t\downarrow\alpha}\omega_{2}(t)=0=\lim_{t\uparrow\beta}\omega_{2}(t).\]
Since $\omega_{2}$ is strictly increasing on $(\alpha,\gamma)$ and
$(\gamma,\beta)$, the point $\gamma$ must really be a pole of $\omega_{2}$
so that $\omega_{2}((\alpha,\gamma))=(0,+\infty)$, $\omega_{2}((\gamma,\beta))=(-\infty,0)$.
We will prove that $\mu_{2}=\delta_{0}$ by showing that $G_{\mu_{2}}$
can be continued analytically across $\mathbb{R}\setminus\{0\}$,
and the continuation is real-valued on $\mathbb{R}\setminus\{0\}$.
Indeed, fix a point $x_{0}\in\mathbb{R}\setminus\{0\}$. There is
a unique $t_{0}\in(\alpha,\beta)$, $t_{0}\neq\gamma$, such that
$\omega_{2}(t_{0})=x_{0}$. Moreover, $\omega_{2}$ is conformal in
a neighborhood of $t_{0}$, so that it has an analytic inverse (with
respect to composition) $\varphi$ defined in a neighborhood $V$
of $x_{0}$, with the property that $\varphi|_{V\cap\mathbb{R}}$
is real-valued. Therefore, we deduce from (2.2) that\[
G_{\mu_{2}}(w)=G_{\mu_{1}\boxplus\mu_{2}}(\varphi(w)),\qquad w\in V\cap\mathbb{C}^{+}.\]
Now, $G_{\mu_{1}\boxplus\mu_{2}}$ is analytic in a neighborhood of
$t_{0}$, and therefore the composition $G_{\mu_{1}\boxplus\mu_{2}}\circ\varphi$
continues analytically to a neighborhood of $x_{0}$. This continuation
is real-valued in an interval around  $x_{0}$ since $G_{\mu_{1}\boxplus\mu_{2}}$
is real-valued in an interval around $t_{0}$.
\end{proof}

\section{Multiplicative Free Convolution on $\mathbb{R}_{+}$}

Given a measure $\mu$ on $\mathbb{R}_{+}=[0,+\infty)$, different
from $\delta_{0}$, we set \begin{equation}
\psi_{\mu}(z)=\int\frac{tz}{1-tz}\, d\mu(t),\label{eq:3.1}\end{equation}
and \begin{equation}
\eta_{\mu}(z)=\frac{\psi_{\mu}(z)}{1+\psi_{\mu}(z)},\qquad z\in\mathbb{C}\setminus\mathbb{R}_{+}.\label{eq:3.2}\end{equation}
The measure $\mu$ is determined by the function $\psi_{\mu}(z)$
since $z(\psi_{\mu}(z)+1)=G_{\mu}\left(\frac{1}{z}\right)$. Note
that the function $\eta_{\mu}$ is characterized by the properties
that $\eta_{\mu}(\bar{z})=\overline{\eta_{\mu}(z)}$, $\lim_{t\uparrow0}\eta_{\mu}(t)=0$,
and $\arg\eta_{\mu}(z)\in[\arg z,\pi)$ for all $z\in\mathbb{C}\setminus\mathbb{R}_{+}$.
For two such measures $\mu_{1},\mu_{2}$, their free multiplicative
convolution $\mu_{1}\boxtimes\mu_{2}$ is characterized by the relation\begin{equation}
\eta_{\mu_{1}\boxtimes\mu_{2}}^{-1}(w)=\frac{1}{w}\eta_{\mu_{1}}^{-1}(w)\eta_{\mu_{2}}^{-1}(w),\qquad w<0.\label{eq:3.3}\end{equation}
As in the case of additive free convolution, the function $\eta_{\mu_{1}\boxtimes\mu_{2}}$
is subordinated to $\eta_{\mu_{j}}$ (see \cite{Biane}). More precisely,
there exist analytic functions $\omega_{1},\omega_{2}:\,\mathbb{C}\setminus\mathbb{R}_{+}\to\mathbb{C}\setminus\mathbb{R}_{+}$
such that\begin{equation}
\eta_{\mu_{1}\boxtimes\mu_{2}}=\eta_{\mu_{1}}\circ\omega_{1}=\eta_{\mu_{2}}\circ\omega_{2},\label{eq:3.4}\end{equation}
 and one can also rewrite (3.3) as\begin{equation}
\eta_{\mu_{1}\boxtimes\mu_{2}}(z)=\frac{1}{z}\omega_{1}(z)\omega_{1}(z),\label{eq:3.5}\end{equation}
 for all $z\in\mathbb{C}\setminus\mathbb{R}_{+}$.  The functions
$\omega_{1},\omega_{2}$ are uniquely determined, and they have the
following properties:

\begin{enumerate}
\item $\lim_{t\uparrow0}\omega_{j}(t)=0,$ $j=1,2$.
\item $\arg z\leq\arg\omega_{j}(z)<\pi$ for all $z\in\mathbb{C}^{+}$,
$j=1,2$. 
\item $\omega_{j}(\bar{z})=\overline{\omega_{j}(z)}$ for $z\in\mathbb{C}\setminus\mathbb{R}_{+}$,
$j=1,2$.
\end{enumerate}
The analogue of Theorem 2.1 for free multiplicative convolution is
proved in \cite{Serban}.

\begin{thm}
Let $\alpha>0$ be an atom of the measure $\mu_{1}\boxtimes\mu_{2}$.
Then
\begin{enumerate}
\item the limits \[
\frac{1}{\alpha_{j}}=\lim_{\varepsilon\downarrow0}\omega_{j}\left(\frac{1}{\alpha}+i\varepsilon\right),\; j=1,2,\]
exist.
\item $\alpha_{1}\alpha_{2}=\alpha$.
\item $\mu_{1}(\{\alpha_{1}\})+\mu_{2}(\{\alpha_{2}\})=(\mu_{1}\boxtimes\mu_{2})(\{\alpha\})+1.$
\item \[
\lim_{\varepsilon\downarrow0}\frac{\omega_{j}\left(\frac{1}{\alpha}+i\varepsilon\right)-\frac{1}{\alpha_{j}}}{i\varepsilon}=\frac{\mu_{j}(\{\alpha_{j}\})}{(\mu_{1}\boxtimes\mu_{2})(\{\alpha\})},\; j=1,2.\]

\end{enumerate}
\end{thm}
Assume now neither $\mu_{1}$ nor $\mu_{2}$ is a point mass, $\alpha,\beta\in(0,+\infty)$
are atoms of $\mu_{1}\boxtimes\mu_{2}$, $\alpha<\beta$, and $(\mu_{1}\boxtimes\mu_{2})((\alpha,\beta))=0$.
Then we can write $\alpha=\alpha_{1}\alpha_{2}$ and $\beta=\beta_{1}\beta_{2}$
by Theorem 3.1. As in the additive case, we may assume that $\alpha_{2}=\beta_{2}=1$
so that $\alpha_{1}=\alpha$ and $\beta_{1}=\beta$ are atoms of the
measure $\mu_{1}$. 

\begin{thm}
Let $\alpha<\beta$ be two positive real numbers such that $\alpha$
and $\beta$ are both atoms for the measure $\mu_{1}\boxtimes\mu_{2}$.
If $(\mu_{1}\boxtimes\mu_{2})((\alpha,\beta))=0$, then either $\mu_{1}$
or $\mu_{2}$ is a point mass. 
\end{thm}
\begin{proof}
With the above notations, we assume that $\alpha_{1}=\alpha$, $\beta_{1}=\beta$,
and $\alpha_{2}=\beta_{2}=1$. The proof proceeds as that of Theorem
2.5. Thus, assuming that $\mu_{1}$ and $\mu_{2}$ are not point masses,
we show
\begin{enumerate}
\item $\eta_{\mu_{1}\boxtimes\mu_{2}}$ continues meromorphically across
$\left(\frac{1}{\beta},\frac{1}{\alpha}\right)$,
\item $\eta_{\mu_{1}\boxtimes\mu_{2}}$ is real-valued on $\left(\frac{1}{\beta},\frac{1}{\alpha}\right)$,
with the exception of a simple pole $\gamma\in\left(\frac{1}{\beta},\frac{1}{\alpha}\right)$,
\item $\omega_{1}$ and $\omega_{2}$ also have continuation properties
in (1) and (2),
\item $\omega_{2}\left(\left(\frac{1}{\beta},\gamma\right)\right)=(1,+\infty)$,
$\omega_{2}\left(\left(\gamma,\frac{1}{\alpha}\right)\right)=(-\infty,1)$,
\item $\eta_{\mu_{2}}$ is real and analytic on $\mathbb{R}\setminus\{1\}$,
hence $\mu_{2}=\delta_{1}$.
\end{enumerate}
The formula defining $\psi_{\mu_{1}\boxtimes\mu_{2}}(z)$ makes sense
for $z\in\left(\frac{1}{\beta},\frac{1}{\alpha}\right)$, so that
(1) and (2) follow immediately from the assumptions on the measure
$\mu_{1}\boxtimes\mu_{2}$. The proof of (3) is analogous to that
of Lemma 2.3. More precisely, we can use the Nevanlinna representation
for functions entering the identity \[
\log\omega_{1}(z)+\log\omega_{2}(z)=\log\eta_{\mu_{1}\boxtimes\mu_{2}}(z)+\log z,\qquad z\in\mathbb{C}^{+},\]
where the princpal value of the logarithm is used. Property (4) then
follows easily from Lemma 2.4, and the fact that $\omega_{2}$ must
be an increasing function on the intervals $\left(\frac{1}{\beta},\gamma\right)$and
$\left(\gamma,\frac{1}{\alpha}\right)$. Finally, (5) follows from
the relation \[
\eta_{\mu_{1}\boxtimes\mu_{2}}(z)=\eta_{\mu_{2}}(\omega_{2}(z))\]
 by locally inverting $\omega_{2}$ around any point in $\mathbb{R}\setminus\{1\}$.
\end{proof}
It should be emphasized that the above result does not hold when $\alpha=0$.
An example is provided by the measures \[
\mu_{1}=\frac{1}{3}\delta_{0}+\frac{2}{3}\delta_{1},\;\mu_{2}=\frac{2}{3}\delta_{1}+\frac{1}{3}\delta_{2}.\]
The free convolution $\mu=\mu_{1}\boxtimes\mu_{2}$ satisfies $\mu(\{0\})=\mu(\{1\})=1/3$,
while $\mu((0,1))=0$. The easiest way to see this is to view $\mu$
as the distribution of the operator $p_{1}(1+p_{2})p_{1}$, where
$p_{1}$ and $p_{2}$ are freely independent selfadjoint projections
in a $W^{*}$-probability space $(\mathcal{A},\tau)$, and $\tau(p_{1})=\tau(p_{2})=1/3$
(We refer to \cite{VDN} for the notions of a $W^{*}$-probability
space and of free independence.)

\section{Free Multiplicative Convolution on $\mathbb{T}$}

For a probability measure $\mu$ on the unit circle $\mathbb{T}$,
the functions $\psi_{\mu}$ and $\eta_{\mu}$ are defined again by
(3.1) and (3.2), but their domain of definition is now the unit disk
$\mathbb{D}=\{ z\in\mathbb{C}:\,\left|z\right|<1\}$. Assume that
\[
\int_{\mathbb{T}}\zeta\, d\mu_{1}(\zeta)\neq0\neq\int_{\mathbb{T}}\zeta\, d\mu_{2}(\zeta).\]
Then the free multiplicative convolution $\mu_{1}\boxtimes\mu_{2}$
is also characterized by (3.3) in a neighborhood of $w=0$. The subordination
functions $\omega_{1}$, $\omega_{2}$ map $\mathbb{D}$ to $\mathbb{D}$,
$\omega_{1}(0)=\omega_{2}(0)=0$, and relations (3.4) and (3.5) are
satisfied in $\mathbb{D}$. Relation (3.5) is satisfied  even when
$\mu_{1}$ or $\mu_{2}$ has first moment equal to zero. 

It is proved in \cite{Serban} that Theorem 3.1 remains valid in this
context. The only changes needed in the statement are that the limits
must be replaced by radial limits. Thus, the formula in part (1) of
Theorem 3.1 becomes\[
\overline{\alpha_{j}}=\frac{1}{\alpha_{j}}=\lim_{r\uparrow1}\omega_{j}(r\overline{\alpha}),\]
while part (4) becomes \[
\lim_{r\uparrow1}\frac{\overline{\alpha_{j}}-\omega_{j}(r\overline{\alpha})}{(1-r)\overline{\alpha}}=\frac{\mu_{j}(\{\alpha_{j}\})}{(\mu_{1}\boxtimes\mu_{2})(\{\alpha\})}.\]

\begin{thm}
Let $\mu_{1}$ and $\mu_{2}$ be two probability measures on $\mathbb{T}$,
and $I\subset\mathbb{T}$ be an open arc with endpoints $\alpha,\beta$.
If $\alpha$ and $\beta$ are atoms of $\mu_{1}\boxtimes\mu_{2}$,
and $(\mu_{1}\boxtimes\mu_{2})(I)=0$, then either $\mu_{1}$ or $\mu_{2}$
is a point mass.
\end{thm}
\begin{proof}
Write $\alpha=\alpha_{1}\alpha_{2}$ and $\beta=\beta_{1}\beta_{2}$,
where $\alpha_{1},\beta_{1}$ are atoms of $\mu_{1}$ and $\alpha_{2},\beta_{2}$
are atoms of $\mu_{2}$. We may assume that $\alpha_{2}=\beta_{2}=1$
so that $\alpha_{1}=\alpha$ and $\beta_{1}=\beta$. As in the earlier
results, we show that
\begin{enumerate}
\item $\eta_{\mu_{1}\boxtimes\mu_{2}}$ continues analytically across $\overline{I}=\{\overline{\zeta}:\,\zeta\in I\}$,
\item $\left|\eta_{\mu_{1}\boxtimes\mu_{2}}(\zeta)\right|=1$ for all $\zeta\in\overline{I}$,
\item $\omega_{1}$ and $\omega_{2}$ also have the continuation properties
stated in (1) and (2),
\item $\omega_{2}\left(\overline{I}\right)=\mathbb{T}\setminus\{1\}$, 
\item $\eta_{\mu_{2}}$ continues analytically across $\mathbb{T}\setminus\{1\}$,
and $\left|\eta_{\mu_{2}}(\zeta)\right|=1$ for $\zeta\in\mathbb{T}\setminus\{1\}$.
Consequently, $\mu_{2}=\delta_{1}$.
\end{enumerate}
The continuation in (1) and (2) is given directly by the formula defining
$\eta_{\mu_{1}\boxtimes\mu_{2}}$. Observe that the zeros of $\eta_{\mu_{1}\boxtimes\mu_{2}}$
have no accumulation points in $\overline{I}$, and therefore the
Blaschke product $B$ corresponding with these zeros is also analytic
across $\overline{I}$. Let us write the decompositions (see \cite[Chapter II]{Garnett})
\[
\eta_{\mu_{1}\boxtimes\mu_{2}}(z)=B(z)\exp\left(\int_{\mathbb{T}}\frac{z+\zeta}{z-\zeta}\, d\sigma(\zeta)\right),\]
\[
\omega_{j}(z)=B_{j}(z)\exp\left(\int_{\mathbb{T}}\frac{z+\zeta}{z-\zeta}\, d\sigma_{j}(\zeta)\right),\qquad j=1,2,\; z\in\mathbb{D},\]
where $\sigma$, $\sigma_{1}$, and $\sigma_{2}$ are finite positive
Borel measures on $\mathbb{T}$. Relation (3.5) implies that \[
zB(z)=B_{1}(z)B_{2}(z),\]
for all $z\in\mathbb{D}$, and that \[
\sigma=\sigma_{1}+\sigma_{2}.\]
Thus, the Blaschke products $B_{1}$ and $B_{2}$ are also analytic
across $\overline{I}$. Moreover, the fact that $\left|\eta_{\mu_{1}\boxtimes\mu_{2}}(\zeta)\right|=1$
for $\zeta\in\overline{I}$ implies that $\sigma(\overline{I})=0$.
We deduce that $\sigma_{1}(\overline{I})=\sigma_{2}(\overline{I})=0$,
and this implies property (3) above. Note that property (3) implies
that $\left|\omega_{2}^{\prime}(\zeta)\right|\geq1$ for all $\zeta\in\overline{I}$. 

The proof of (4) follows from the fact that $\lim_{\zeta\in\overline{I},\:\zeta\rightarrow\overline{\alpha}}\omega_{2}(\zeta)=\lim_{\zeta\in\overline{I},\:\zeta\rightarrow\overline{\beta}}\omega_{2}(\zeta)=1$.
To see this, one must use the analogue of Lemma 2.4, which can be
proved by using the conformal equivalence between $\mathbb{C}^{+}$
and $\mathbb{D}$. Finally, (5) follows from the identity \[
\eta_{\mu_{1}\boxtimes\mu_{2}}(z)=\eta_{\mu_{2}}(\omega_{2}(z))\]
by locally inverting $\omega_{2}$ around any point in $\mathbb{T}\setminus\{1\}$. 
\end{proof}

\end{document}